\documentclass[10pt,journal,compsoc]{IEEEtran}

\usepackage[nocompress]{cite}

\usepackage[pdftex]{graphicx}

\usepackage{amsmath}
\usepackage{amsfonts}
\usepackage{amsthm}
\newtheorem{theorem}{Theorem}
\newtheorem*{proof*}{Proof}
\newtheorem{cor}{Corollary}
\newtheorem{lemma}{Lemma}

\usepackage{algorithm,algorithmic}

\usepackage{array}

\usepackage[caption=false,font=footnotesize,labelfont=sf,textfont=sf]{subfig}

\usepackage{stfloats}

\usepackage{url}

\begin{document}

\title{Generating Graphs with Symmetry}

\author{Isaac~Klickstein~and~Francesco~Sorrentino%
\IEEEcompsocitemizethanks{\IEEEcompsocthanksitem I. Klickstein and F. Sorrentino are with the Department
of Mechanical Engineering, University of New Mexico, Albuquerque,
NM, 87131.\protect\\
E-mail: iklick@unm.edu}%
}

\markboth{IEEE Transactions on Network Science and Engineering}%
{Klickstein~and~Sorrentino}

\IEEEtitleabstractindextext{%
\begin{abstract}
In the field of complex networks and graph theory, new results are typically tested on graphs generated by a variety of algorithms such as the Erd\H{o}s-R\'{e}nyi model or the Barab\'{a}si-Albert model.
Unfortunately, most graph generating algorithms do not typically create graphs with symmetries, which have been shown to have an important role on the network dynamics.
Here, we present an algorithm to generate graphs with prescribed symmetries.
The algorithm can also be used to generate graphs with a prescribed equitable partition but possibly without any symmetry.
We also use our graph generator to examine the recently raised question about the relation between the orbits of the automorphism group and a graph's minimal equitable partition.
\end{abstract}

}

\maketitle

\IEEEdisplaynontitleabstractindextext

\IEEEpeerreviewmaketitle

\IEEEraisesectionheading{\section{Introduction}\label{sec:introduction}}

\IEEEPARstart{D}{ue} to the interest in the field of studies on complex networks, a number of network generating algorithms have been proposed. Among these are the Erd\H{o}s-R\'{e}nyi random graph model \cite{erdos1959random},  the Watts-Strogatz small world model \cite{watts1998collective}, the Barab\'{a}si-Albert model \cite{barabasi1999emergence} which generates scale free networks, the static model \cite{goh2001universal} and the configuration model \cite{molloy1995critical} (together with its uncorrelated version \cite{catanzaro2005generation}) which have been used to reproduce scale free networks with given power-law degree distribution exponents, and a number of models that generate networks with assigned degree distribution and degree correlation \cite{newman2002assortative,newman2003mixing}.
However, currently available network generating algorithms very rarely reproduce symmetries.
It has been shown that symmetries are present and play an important role in dynamical systems with topology described by a graph \cite{macarthur2008symmetry,belykh2001cluster,sorrentino2016complete}.
It is thus important to introduce a simple network generating algorithm that can create networks with a desired number of symmetries.\\
\indent
Here, we present an algorithm that can generate graphs with prescribed symmetries.
Additionally, the algorithm is extended to the case one wants to generate random graphs with a prescribed equitable partition \cite{belykh2008cluster}.
Using the fact we are able to generate graphs with prescribed equitable partitions but possibly without symmetries, we can further investigate the recently raised question \cite{siddique2018symmetry} concerning when the minimal equitable partition and the graph symmetries align, and when they do not.\\
\indent 
In section \ref{sec:background} we present definitions of equitable partitions and the orbits of the automorphism group, as well as how we may compress a graph with an equitable partition to its quotient graph.
In section 3 we derive the algorithm and present an example.
Finally, in section 4 we use the algorithm to compare the minimal equitable partition with the partition induced by the orbits of the automorphism group and make concluding remarks in section 5.

\section{Preliminaries}\label{sec:background}
Let $\mathcal{G} = (\mathcal{V},\mathcal{E})$ denote a simple, undirected, unweighted graph with the set of $n$ vertices $\mathcal{V}$ and the set of edges $\mathcal{E} \subset \mathcal{V} \times \mathcal{V}$.
As $\mathcal{G}$ is undirected, each edge is an unordered pair of vertices $(v_i,v_j) = (v_j,v_i)$.
Also, as $\mathcal{G}$ is simple, we do not allow for any self-loops, $(v_i,v_i) \notin \mathcal{E}$, and no multi-edges, i.e., an edge $(v_i,v_j) \in \mathcal{E}$ may only appear once.
A graph $\mathcal{G}$ can be represented as an $n \times n$ adjacency matrix $G$ where $G_{ij} = 1$ if $(v_j,v_i) \in \mathcal{E}$ and $G_{ij} = 0$ otherwise.
As we assume $\mathcal{G}$ is undirected, the adjacency matrix $G$ is symmetric.
A partition $\mathcal{C}$ of the vertices $\mathcal{V}$ satisfies the properties,
\begin{equation}\label{eq:partition}
  \mathcal{C} = \left\{\mathcal{C}_i \subset \mathcal{V} \left| \mathcal{C}_i \cap \mathcal{C}_j = \emptyset, \quad |\mathcal{C}_i| = n_i, \quad \sum_{i=1}^{p} n_i = n \right.\right\}
\end{equation}
where we call $\mathcal{C}_i$ the $i$th cluster, $i = 1, \ldots, p$, $n_i$ is the number of nodes in cluster $\mathcal{C}_i$, and $n$ is the total number of nodes in the graph. 
An equitable partition (or balanced coloring \cite{belykh2011mesoscale}) $\mathcal{C}$ of a graph $\mathcal{G}$ is a partition with the additional property that,
\begin{equation}\label{eq:equi}
  \sum_{v_a \in \mathcal{C}_k} G_{ia} = \sum_{v_a \in \mathcal{C}_k} G_{ja}, \quad \begin{aligned}
    &\forall v_i,v_j \in \mathcal{C}_{\ell}\\ &\forall \mathcal{C}_k, \mathcal{C}_{\ell} \in \mathcal{C}
  \end{aligned}
\end{equation}
The relation in Eq. \eqref{eq:equi} states that if two vertices are in the same cluster, $v_i,v_j \in \mathcal{C}_{\ell}$, then they must be adjacent to the same number of vertices in each of the clusters.
There are two balanced colorings we are particularly interested in, the minimal balanced coloring (MBC) and the orbits of the automorphism group (OAG) of the graph $\mathcal{G}$.
The MBC \cite{belykh2011mesoscale} is the balanced coloring $\mathcal{C}$ of a graph $\mathcal{G}$ that solves the optimization problem,
\begin{equation}
  \begin{aligned}
    \min && &p\\
    \text{s.t.} && &|\mathcal{C}| = p\\
    && &\mathcal{C} \text{ is a balanced coloring}
  \end{aligned}
\end{equation}
The OAG is best defined using the symmetry group of permutations of the graph $\mathcal{G}$.
A permutation of the vertices $\mathcal{V}$ is a bijection $\pi : \mathcal{V} \mapsto \mathcal{V}$ which can be thought of as a shuffling of the vertices, i.e., no vertices are removed or created.
Each permutation can be represented by an $n\times n$ permutation matrix $P$ where $P_{ij} = 1$ if $\pi(v_j) = v_i$ and $P_{ij} = 0$ otherwise.
Note that by this definition, $P$ is orthonormal so that $PP^T = I$, the appropriately dimensioned identity matrix.
A generic permutation of the vertices in a graph will alter the set of edges $\mathcal{E}$.
A symmetry of a graph is a permutation such that $\mathcal{E}$ remains unchanged after the permutation.
Algebraically, a symmetric permutation must be $G$ invariant.
\begin{equation}
  PGP^T = G
\end{equation}
The set of these symmetric permutations, or simply symmetries, is a group under permutation composition called the automorphism group, $\text{Aut}(\mathcal{G})$ \cite{lauri2016topics}.
The automorphism group induces an equitable partition of the nodes in the graph called the orbits of the automorphism group (OAG).
In general, the equitable partitions OAG and MBC are not equal \cite{siddique2018symmetry,schaub2016graph,kudose2009equitable}.
Coupled dynamical systems whose underlying graph has a non-trivial OAG partition can exhibit complex behavior where the nodes in the same cluster may behave similarly even if they are not directly connected \cite{nicosia2013remote,pecora2014cluster}.\\
\indent
The following theorem states the requirement for $\mathcal{C} = \mathcal{O}$.
\begin{theorem}
  Let $\phi_{ij}$ denote the operation of swapping two nodes $v_i$ and $v_j$ such that $v_i,v_j \in \mathcal{C}_k$.
  Also, let $\mathcal{O}$ denote the partition of the nodes of graph $\mathcal{G} = (\mathcal{V},\mathcal{E})$ induced by its automorphism group $\text{Aut}(\mathcal{G})$ and let $\mathcal{C}$ denote the partition of the nodes of the same graph induced by its MBC.
  Then $\mathcal{C} = \mathcal{O}$ if and only if for every $\phi_{ij}$, one can construct a permutation
  \begin{equation}
    \pi = \phi_{ij} \phi_{k_1k_1'} \phi_{k_2k_2'} \ldots \phi_{k_sk_s'}
  \end{equation}
  such that $k_{\ell} \neq i,j$ and $k_{\ell}' \neq i,j$ for $1 \leq \ell \leq s$ and $\pi \in \text{Aut}(\mathcal{G})$.
\end{theorem}
The theorem holds by the definition of the automorphism group, and is useful as a tool to check whether or not one should expect the minimal balanced coloring and the orbits of the automorphism group to coincide.\\
\indent
Typically, large random graphs generated with the Erd\H{o}s-Renyi model, the Watts-Strogatz model, the Barab{\'a}si-Albert model, and most others will not have non-trivial equitable partitions, that is, $|\mathcal{C}| =|\mathcal{V}|$.
If one generates a random network using any of these methods one will not see the effect that equitable partitions can have on the system dynamics.
However, real networks are often characterized by a large number of symmetries \cite{macarthur2008symmetry}.
This prompts us to study in this paper a procedure to generate large graphs with an assigned number of symmetries.\\
\begin{figure}
  \centering
  \includegraphics[width=2.5in]{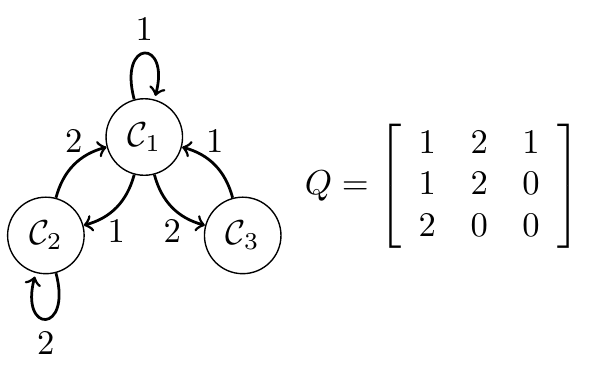}
  \caption{An example of a quotient graph with three vertices (cluster) $\mathcal{C} = \{\mathcal{C}_1,\mathcal{C}_2,\mathcal{C}_3\}$.
  The self-loops are $\mathcal{D}(\mathcal{C}_1) = 1$, $\mathcal{D}(\mathcal{C}_2) = 2$, and $\mathcal{D}(\mathcal{C}_3) = 0$.
  There are two edges, $\mathcal{F}_1 = (\mathcal{C}_1,\mathcal{C}_2)$ and $\mathcal{F}_2 = (\mathcal{C}_1,\mathcal{C}_3)$.
  The first edge has weights $\mathcal{W}(\mathcal{F}_1,0) = 1$ and $\mathcal{W}(\mathcal{F}_1,1) = 2$, i.e., from $\mathcal{C}_2$ to $\mathcal{C}_1$ and from $\mathcal{C}_1$ to $\mathcal{C}_2$, respectively.
  The second edge has weights $\mathcal{W}(\mathcal{F}_2,0) = 2$ and $\mathcal{W}(\mathcal{F}_2,1) = 1$.
  Also shown is the adjacency matrix of this quotient graph, $Q$.}
  \label{fig:quotient}
\end{figure}
An equitable partition $\mathcal{C}$ of a graph $\mathcal{G}$ can be represented as a quotient graph $\mathcal{Q} = (\mathcal{C},\mathcal{D},\mathcal{F},\mathcal{W})$.
An example of a quotient graph with a discussion of its four components is shown in Fig. \ref{fig:quotient}.
Each vertex in the quotient graph $\mathcal{C}_k \in \mathcal{C}$, $k=1,\ldots,p$ represents the set of vertices in the same cluster in the original graph $\mathcal{G}$.
The self-loop magnitudes $\mathcal{D} : \mathcal{C} \mapsto \mathbb{N}$ is the number of edges each node $v_i \in \mathcal{C}_k$ receives from the other nodes in $\mathcal{C}_k$.
\begin{equation}
  \mathcal{D}(\mathcal{C}_k) = \sum_{v_j \in \mathcal{C}_k} G_{ij}, \quad v_i \in \mathcal{C}_k
\end{equation}
Note that by the definition of an equitable partition, the particular choice of $v_i \in \mathcal{C}_k$ does not affect the value of $\mathcal{D}(\mathcal{C}_k)$.
The edges $\mathcal{F} \subset \mathcal{C} \times \mathcal{C}$ represent those pairs of clusters with edges passing between them.
Each edge has two weights $\mathcal{W} : \mathcal{F} \times \{0,1\} \mapsto \mathbb{N}^+$ 
\begin{equation}
  \begin{aligned}
    \mathcal{W}((\mathcal{C}_k,\mathcal{C}_{\ell}),0) = \text{ weight from $\ell$ to $k$}, \quad k < \ell\\
    \mathcal{W}((\mathcal{C}_k,\mathcal{C}_{\ell}),1) = \text{ weight from $k$ to $\ell$}, \quad k < \ell
  \end{aligned}
\end{equation}
In words, for each edge $\mathcal{F}_i = (\mathcal{C}_k,\mathcal{C}_{\ell})$, $\mathcal{W}(\mathcal{F}_i,0)$ is the weight to the lower indexed vertex from the higher indexed vertex and $\mathcal{W}(\mathcal{F}_i,1)$ is the weight to the higher indexed vertex from the lower indexed vertex in the quotient graph.\\
\indent
The quotient graph $\mathcal{Q}$ can be represented as a $p \times p$ matrix $Q$ with entries,
\begin{equation}
  Q_{ij} = \left\{ \begin{aligned}
    &\mathcal{D}(\mathcal{C}_i), && i = j\\
    &\mathcal{W}((\mathcal{C}_i,\mathcal{C}_j),0), && (\mathcal{C}_i,\mathcal{C}_j) \in \mathcal{F}, && i < j\\
    &\mathcal{W}((\mathcal{C}_i,\mathcal{C}_j),1), && (\mathcal{C}_i,\mathcal{C}_j) \in \mathcal{F}, && i > j\\
    &0, && \text{otherwise}
  \end{aligned} \right.
\end{equation}
One should read $Q_{ij}$ as `vertices in cluster $\mathcal{C}_i$ receive $Q_{ij}$ edges from vertices in cluster $\mathcal{C}_j$' and $Q_{ii}$ as `each vertex in cluster $\mathcal{C}_i$ receives $Q_{ii}$ edges from other vertices in $\mathcal{C}_i$'.
We remark that in this paper, we always distinguish between self-loops $Q_{ii}$ (connections from a vertex to itself in the quotient graph) and edges $Q_{ij}$ (connections between two different vertices in the quotient graph), as they will be treated very differently in the forthcoming derivations.\\
\indent
The approach of this paper is to generate a full network given knowledge of its quotient graph.
As we will see, there are certain prerequisites a quotient graph must satisfy in order to be \emph{feasible}, i.e., for the existence of a transformation mapping the quotient graph to a corresponding full graph.
The quotient graph $\mathcal{Q}$ is unique for a given graph $\mathcal{G}$ with equitable partition $\mathcal{C}$.
On the other hand, a single, feasible, quotient graph represents an infinite number of original graphs.
Thus, reconstructing a graph $\mathcal{G}$ from a feasible quotient graph $\mathcal{Q}$ is not unique and we must select from this infinite set.\\
\indent 
Our procedure illustrated in what follows is based on three steps: (i) we select a quotient graph and assess whether it is feasible, (ii) from a feasible quotient graph, we determine an equitable partition of the network nodes and (iii), we wire the edges of the network so as to ensure that the MBC and OAG coincide.
Each one of these three steps is presented in detail in sections \ref{ssec:feas}, \ref{ssec:ilp}, and \ref{ssec:wire}, respectively.

\section{Results}\label{sec:results}
Before defining the algorithm, we present some useful results from the graph theory literature.
\begin{theorem}[Erd\H{o}s-Gallai \cite{choudum1986simple}]\label{thm:EG}
  Let $a = (a_1,a_2,\ldots,a_n)$ be a non-increasing sequence of non-negative integers.
  The sequence $a$ is realizable as the degree sequence of an undirected simple graph (i.e., one with no self-loops or multi-edges) if and only if,
  \begin{enumerate}
    \item $\sum_{i=1}^n a_i$ is even, and
    \item for $1 \leq k \leq n$,
    \begin{equation}
      \sum_{i=1}^k a_i \leq k(k-1) + \sum_{i=k+1}^n \min \{a_i,k\}
    \end{equation}
  \end{enumerate}
\end{theorem}
\begin{cor}\label{cor:EG}
  If the degree sequence consists of a constant, $a = (r,r,\ldots,r)$, then the two conditions in Thm. \ref{thm:EG} can be written in terms of the length of the sequence $n$.
  \begin{enumerate}
    \item If $r$ is even, then $n$ may be even or odd. If $r$ is odd, then $n$ must be even.
    \item From the case $r \leq k$ in the second condition, it can be shown $n \geq r+1$. For the case $r > k$, the second condition is trivially satisfied.
  \end{enumerate}
\end{cor}
\begin{theorem}[Gale-Ryser \cite{gale1957theorem}]\label{thm:GR}
  Let $a = (a_1,a_2,\ldots,a_{n_1})$ and $b = (b_1,b_2,\ldots,b_{n_2})$ be two non-increasing sequences of non-negative integers.
  The sequences $a$ and $b$ can be realized as the degree sequences of a simple bipartite graph if and only if,
  \begin{equation}
    \sum_{i=1}^k a_i \leq \sum_{i=1}^{n_2} \min\{b_i,k\}, \quad 1 \leq k \leq n_1
  \end{equation}
  or, equivalently,
  \begin{equation}
    \sum_{i=1}^k b_i \leq \sum_{i=1}^{n_1} \min \{a_i,k\}, \quad 1 \leq k \leq n_2
  \end{equation}
\end{theorem}
\begin{cor}\label{cor:GR}
  If $a = (r_1,r_1,\ldots,r_1)$ and $b = (r_2,r_2,\ldots,r_2)$ are two sequences of a constant integers $r_1$ and $r_2$ of lengths $n_1$ and $n_2$, respectively, then the conditions in Thm. \ref{thm:GR} can be rewritten as,
  \begin{equation}
    r_1 \leq n_2, \quad r_2 \leq n_1
  \end{equation}
  and
  \begin{equation}
    r_1n_1 = r_2n_2
  \end{equation}
\end{cor}
Also useful will be the following result.
\begin{lemma}\label{lem:pos}
  Let $A \in \mathbb{R}^{m \times n}$ be a matrix such that each row consists of a single positive entry, a single negative entry, and the remaining entries are all zero, and no columns of $A$ consist entirely of zeros.
  If $\mathcal{N}(A) \neq \emptyset$, there always exists a positive vector $\textbf{x}$, that is a vector with all strictly positive entries,  such that $\textbf{x} \in \mathcal{N}(A)$.
\end{lemma}
\begin{proof*}
  By assumption, let there exist at least one solution (besides $\textbf{x} = \boldsymbol{0}$) of $A \textbf{x} = \boldsymbol{0}$  and let $\textbf{a}_i \in \mathbb{R}^n$ be the $i$th row of $A$ with $a_{i,j} > 0$ and $\textbf{a}_{i,k}<0$ by the construction of $A$.
  First, assume the entries of $\textbf{x}$ are $x_j > 0$ and $x_k < 0$.
  Then, $\textbf{a}_i^T \textbf{x} = \textbf{a}_{i,j}x_j + \textbf{a}_{i,k} x_k > 0$.
  On the other hand, if $x_j < 0$ and $x_k > 0$, then $\textbf{a}_i^T \textbf{x} < 0$.
  By contradiction, if $\textbf{x} \in \mathcal{N}(A)$, then either $x_j > 0$ and $x_k > 0$ or $x_j < 0$ and $x_k < 0$.
  If we find a negative vector $\textbf{x}$, that is a vector with all strictly negative entries, such that $\textbf{x} \in \mathcal{N}(A)$, then obviously $-\textbf{x} \in \mathcal{N}(A)$ as well, and $-\textbf{x}$ is a positive vector.
\end{proof*}
\begin{lemma}\label{lem:int}
  Let $A$ be defined as in Lemma \ref{lem:pos} with the additional constraint that all non-zero entries in $A$ are rational numbers.
  If $\mathcal{N}(A) \neq \emptyset$, then there exists a positive integer solution $\textbf{x} \in \mathbb{Z}^n$ such that $\textbf{x} \in \mathcal{N}(A)$.
\end{lemma}
\begin{proof*}
  By assumption, there exists some vector $\textbf{y} \in \mathbb{R}^n$ such that $A \textbf{y} = \boldsymbol{0}$.
  We can solve for a positive vector $\textbf{y}$ by Gaussian elimination which involves only elementary operations so that each entry in $y_i = p_i / q_i$ where $p_i$ and $q_i$ are integers, i.e., $y_i$ is a rational number.
  Define $k = \prod_{i=1}^n q_i$ so that $ky_i$ is an integer and thus $\textbf{x} = k \textbf{y}$ is a positive integer vector.
  Clearly then, $k A \textbf{y} = A (k\textbf{y}) = A \textbf{x} = \boldsymbol{0}$ so that $\textbf{x} \in \mathcal{N}(A)$.
\end{proof*}

As we will see next, given a quotient graph $\mathcal{Q}$ from which we are to construct a symmetric unweighted graph $\mathcal{G}$ we must perform three tasks; (i) verify that $\mathcal{Q}$ is a feasible quotient graph, and if it is (ii) determine the cardinality of each $|\mathcal{C}_i| = n_i$, $i= 1,\ldots,p$, after which, finally, (iii) we wire the edges according to $\mathcal{F}$ and $\mathcal{W}$.
\subsection{Feasibility of the Quotient Graph}\label{ssec:feas}
Using Corollaries \ref{cor:EG} and \ref{cor:GR} we can construct the set of requirements in terms of the cardinalities $n_i$ for $\mathcal{Q}$ to be a feasible quotient graph.
\begin{subequations}
\begin{align}
    \label{eq:lb}
    n_i &\geq \max \left\{ \max\limits_{(\mathcal{C}_j,\mathcal{C}_i) \in \mathcal{F}} Q_{ji}, Q_{ii}+1 \right\}\\
    \label{eq:even}
    \bmod (Q_{ii} n_i, 2) &= 0\\
    \label{eq:cons}
    Q_{ij} n_i &= Q_{ji} n_j, \quad \forall (\mathcal{C}_i,\mathcal{C}_j) \in \mathcal{F}
\end{align}
\end{subequations}
To enforce the constraint in Eq. \eqref{eq:even}, we define a new set of variables $x_i$ for each cluster $\mathcal{C}_i$.
\begin{equation}\label{eq:x2n}
  \begin{aligned}
    x_i = n_i, \quad Q_{ii} \bmod 2 = 0\\
    2x_i = n_i, \quad Q_{ii} \bmod 2 = 1
  \end{aligned}
\end{equation}
The definition of $x_i$ modifies the constraint in Eq. \eqref{eq:lb} slightly.
\begin{equation}\label{eq:xlb}
  \begin{aligned}
    x_i \geq \left(1 - \frac{Q_{ii} \bmod 2}{2}\right) \max \left\{ \max\limits_{(\mathcal{C}_j,\mathcal{C}_i) \in \mathcal{F}} Q_{ji}, Q_{ii}+1 \right\}
  \end{aligned}
\end{equation}
For notational ease, let the lower bound of $x_i$ in Eq. \eqref{eq:xlb} be defined as $x^L_i$.
The set of constraints in Eq. \eqref{eq:cons} can be combined into a system of linear equations $A \textbf{x} = \boldsymbol{0}$ where $A$ is a $|\mathcal{F}| \times |\mathcal{C}|$ matrix with each row corresponding to an edge $\mathcal{F}_{i} = (\mathcal{C}_j,\mathcal{C}_k)$ with entries
\begin{equation}\label{eq:Adef}
  A_{ij} = \left\{ \begin{aligned}
    Q_{jk}, && &\text{if} \quad j < k \quad \text{and} \quad Q_{jj} \bmod 2 = 0\\
    2Q_{jk}, && &\text{if} \quad j < k \quad \text{and} \quad Q_{jj} \bmod 2 = 1\\
    -Q{kj}, && &\text{if} \quad j > k \quad \text{and} \quad Q_{jj} \bmod 2 = 0\\
    -2Q_{kj}, && &\text{if} \quad j > k \quad \text{and} \quad Q_{jj} \bmod 2 = 1\\
    0, && &\text{otherwise}
  \end{aligned} \right.
\end{equation}
The matrix $A$ defined in Eq. \eqref{eq:Adef} is of the form presented in Lemma \ref{lem:pos}.
By Lemmas \ref{lem:pos} and \ref{lem:int}, if $\mathcal{N}(A) \neq \emptyset$, then there exists an integer solution $\textbf{x}$ such that $A \textbf{x} = \boldsymbol{0}$, $\textbf{x} \geq \textbf{x}^L$ and we can then reconstruct $\textbf{n}$ using Eq. \eqref{eq:x2n}.
Note that if we find one solution $\textbf{n}$, then any integer multiple $s\textbf{n}$, $s \geq 1$, is also a solution.
On the other hand, if $\mathcal{N}(A) = \emptyset$, then there is no solution $\textbf{n}$, which leads to our first main result.
\begin{theorem}
  The quotient graph $\mathcal{Q} = (\mathcal{C},\mathcal{D},\mathcal{F},\mathcal{W})$ is a feasible quotient graph if and only if the matrix $A$ defined in Eq. \eqref{eq:Adef} has a non-trivial null space.
\end{theorem}
This result provides a framework with which one may construct feasible quotient graphs.
One first freely chooses $p = |\mathcal{C}|$, the number of clusters, and arbitrarily assigns $Q_{ii}$, $i = 1,\ldots,p$.
One may also freely add $p-1$ edges $\mathcal{F}_i = (\mathcal{C}_j,\mathcal{C}_k)$ with weights $\mathcal{W}(\mathcal{F}_i,0))$ and $\mathcal{W}(\mathcal{F}_i,1)$.
At this point, the maximum possible rank of $A$ is $(p-1)$, and so the dimension of the null space is at least one.
Any additional edges may be added only if they are linear combinations of the first $(p-1)$ rows of $A$ and satisfy the condition that each row of $A$ has one positive integer and one negative integer with all other entries zero.
%
\subsection{The Equitable Partition of the Nodes}\label{ssec:ilp}
After verifying that $\mathcal{Q}$ is a feasible quotient graph simply by checking that $\mathcal{N}(A) \neq \emptyset$, the question remains how one may find $\textbf{n}$, the cardinalities of each cluster $n_1,n_2,\ldots,n_{p}$.
One option is to compute a basis for the null space of $A$ and then proceed as in the proof of Lemma \ref{lem:int}.
This procedure may lead to very large $n_i$ during the scaling of the rational solution $\textbf{y}$ to an integer solution.
We propose an alternative approach by formulating an integer linear program (ILP) that can yield a minimal realization by which we mean the resulting graph $\mathcal{G}$ has the fewest number of nodes,
\begin{equation}\label{eq:ilp}
  \begin{aligned}
    \min && &\textbf{c}^T \textbf{x}\\
    \text{s.t.} && &A \textbf{x} = \boldsymbol{0}\\
    && &\textbf{x} \geq \textbf{x}^L
  \end{aligned}
\end{equation}
The vector $\textbf{c}$ is any strictly positive vector whose particular choice will select different integer vectors $\textbf{x}$ in the null space of $A$.

By ensuring $\mathcal{Q}$ is a feasible quotient graph, we can be certain the dimension of the null space of $A$ is at least one, and Eq. \eqref{eq:ilp} has a solution.
To solve instances of Eq. \eqref{eq:ilp} we use the branch and cut algorithm implemented in COIN-OR's package CBC \cite{lougee2003common}.
Once we have solved Eq. \eqref{eq:ilp}, $\textbf{x}^*$, we apply Eq. \eqref{eq:x2n} to compute $\textbf{n}$, and we may scale it larger by any integer multiple $s$ if one is interested in a graph with a larger number of vertices.\\
\subsection{Wiring the Edges}\label{ssec:wire}
Wiring the edges must be done in such a way to ensure that the desired number of edges between clusters dictated by the quotient graph is satisfied and the maximum number of symmetries is present so the OAG and the MBC coincide.
To simplify the notation in what follows, for each node $v_j \in \mathcal{V}$, $j = 1, \ldots, n$, in the graph, we assign a cluster specific label, $u_i \in \mathcal{C}_k$, $i = 0, \ldots, n_k-1$.
Each self-loop and each edge in the quotient graph can be handled separately.
First we will discuss the wiring to satisfy each self-loop, and second we will discuss the wiring to satisfy each edge in the quotient graph.
\subsubsection{Intra-Cluster Edges}
Let $u_i \in \mathcal{C}_k$, $i = 0,\ldots,n_k-1$ be a node in cluster $k$ and assume $Q_{kk} > 0$ (otherwise if $Q_{kk} = 0$ there are no intra-cluster edges to add for $\mathcal{C}_k$).
In what follows, we must distinguish between the two cases that $Q_{kk}$ is even or odd.
If $Q_{kk}$ is even, we add the edges,
\begin{equation}\label{eq:intraeven}
  \left\{ \begin{aligned}
    &\left( u_i, u_{(i+j) \bmod n_k} \right), && j = 1,\ldots,Q_{kk}/2\\
    &\left( u_i, u_{(i+n_k-j) \bmod n_k} \right), && j = 1,\ldots,Q_{kk}/2
  \end{aligned} \right.
\end{equation}
From the formulation of the ILP, $n_k > Q_{kk}$ so that each node $u_i$ is connected to $Q_{kk}$ other nodes in the same cluster.
We also show that the set of edges in Eq. \eqref{eq:intraeven} is symmetric, that is, if $(u_i,u_{i'})$ is an edge, then $(u_{i'},u_i)$ is an edge as well.
From the first line of Eq. \eqref{eq:intraeven} we see that $i' = (i+j) \bmod n_k$.
Substituting this expression for $i$ into the second line of Eq. \eqref{eq:intraeven}, we must show that,
\begin{equation}
  \begin{aligned}
    &\left[ \left( i + j\right)\bmod n_k + n_k - j \right] \bmod n_k\\
    &= \left[ i + j - qn_k + n_k - j \right] \bmod n_k\\
    &= \left[ i + (1-q) n_k \right] \bmod n_k = i \bmod n_k = i
  \end{aligned}
\end{equation}
In a similar fashion one can show that if $i' = (i + n_k - j) \bmod n_k$ from the second line of Eq. \eqref{eq:intraeven}, we can substitute $i'$ for $i$ in the first line of Eq. \eqref{eq:intraeven}.
Thus, for $0 \leq i \leq n_k-1$, there are $n_k Q_{kk}/2$ unique edges created by Eq. \eqref{eq:intraeven} as every pair appears twice.
If $Q_{kk}$ is odd, we add the following edges,
\begin{equation}\label{eq:intraodd}
  \left\{ \begin{aligned}
    &\left( u_i, u_{(i+j)\bmod n_k} \right), && j = 1, \ldots, (Q_{kk}-1)/2\\
    &\left( u_i, u_{(i+n_k-j) \bmod n_k} \right), && j = 1, \ldots, (Q_{kk}-1)/2\\
    &\left( u_i, u_{(i + n_k/2) \bmod n_k} \right)
  \end{aligned} \right. 
\end{equation}
By the formulation of the ILP, if $Q_{kk}$ is odd, $n_k$ is even so $n_k/2$ is an integer.
Showing that the set of edges in Eq. \eqref{eq:intraodd} leads to each node $u_i \in \mathcal{C}_k$ having intra-cluster degree $Q_{kk}$ and that the set of edges is symmetric is very similar to the proof of Eq. \eqref{eq:intraeven} and thus we do not include it here.
\begin{figure*}
  \centering
  \includegraphics[width=7.3in]{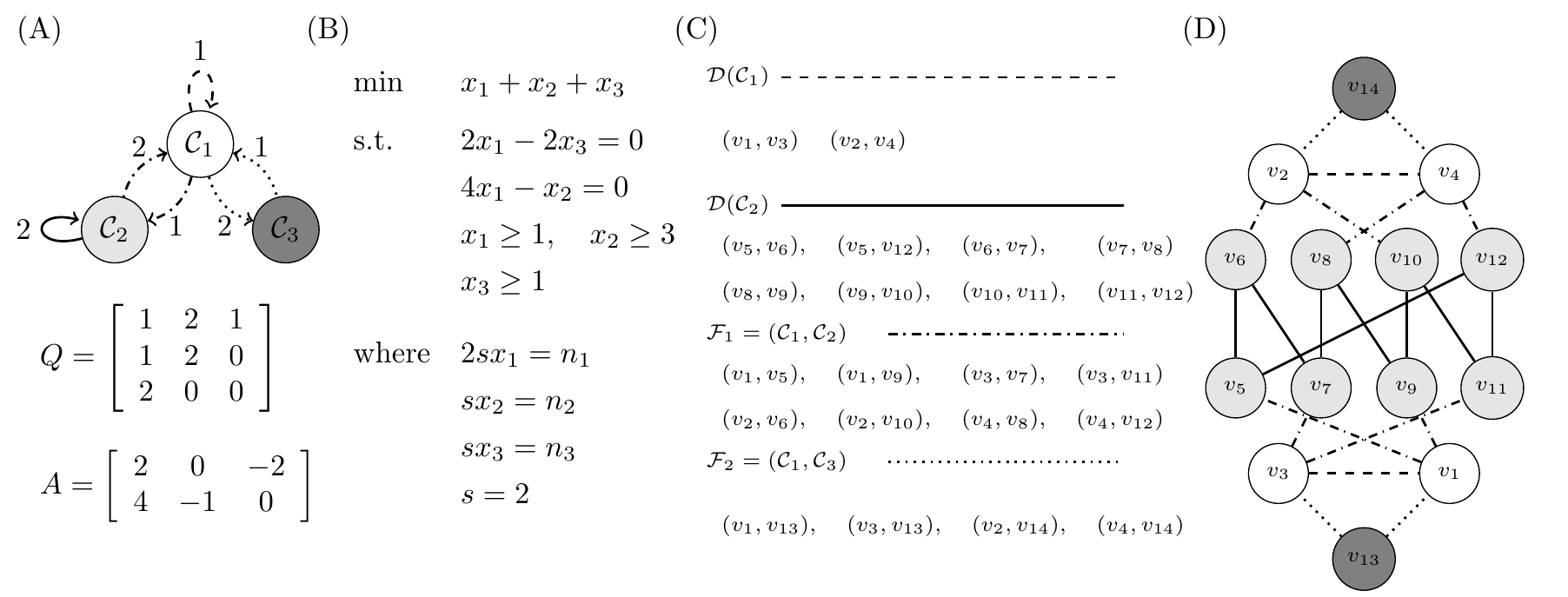}
  \caption{An example of the process from quotient graph to full graph.
  (A) A graphical depiction of the quotient graph $\mathcal{Q}$ with three vertices $\mathcal{C}_1$, $\mathcal{C}_2$, and $\mathcal{C}_3$.
  Vertices $\mathcal{C}_1$ and $\mathcal{C}_2$ have self-loops $\mathcal{D}(\mathcal{C}_1) = 1$ and $\mathcal{D}(\mathcal{C}_2) = 2$.
  There are two edges, $\mathcal{F}_1 = (\mathcal{C}_1,\mathcal{C}_2)$ and $\mathcal{F}_2 = (\mathcal{C}_1,\mathcal{C}_3)$.
  The first edge has weights $\mathcal{W}(\mathcal{F}_1,0) = 2$ and $\mathcal{W}(\mathcal{F}_1,1) = 1$ while the second edge has weights $\mathcal{W}(\mathcal{F}_2,0) = 1$ and $\mathcal{W}(\mathcal{F}_2,1) = 2$.
  Also shown is the weighted quotient adjacency matrix $Q$ and the matrix $A$ that appears as the linear constraints in the ILP in Eq. \eqref{eq:ilp}.
  Note that the dimension of the null space of $A$ is one so that there exists a solution to $A \textbf{x} = \boldsymbol{0}$.
  (B) The ILP that must be solved to find $x_1$, $x_2$, and $x_3$.
  Note that we will scale the population by $s = 2$.
  As $\mathcal{D}(\mathcal{C}_1)$ is odd, we must set $2sx_1 = n_1$ while the other populations are $sx_2 = n_2$ and $sx_3 = n_3$.
  The consistent populations found are $n_1 = 4$, $n_2 = 8$, and $n_3 = 2$.
  (C) The sets of edges corresponding to each quotient graph self-loop and edge are computed.
  The edges prescribed by the self-loop in $\mathcal{C}_1$ are found using Eq. \eqref{eq:intraodd}.
  The edges prescribed by the self-loop in $\mathcal{C}_2$ are found using Eq. \eqref{eq:intraeven}.
  The edges in the full graph prescribed by the two edges in the quotient graph $\mathcal{F}_1$ and $\mathcal{F}_2$ are found using Eq. \eqref{eq:uw}.
  (D) A diagram of the resulting graph consisting of $n=14$ vertices. The vertices are shaded according to their cluster and the edges are textured according to their corresponding edge in the quotient graph shown in the diagram in (A).
   }
  \label{fig:example}
\end{figure*}
\subsubsection{Inter-Cluster Edges}
Let $u_i \in \mathcal{C}_k$, $i = 0,\ldots,n_k-1$ and $w_j \in \mathcal{C}_{\ell}$, $j = 0 \ldots, n_{\ell}-1$ be nodes in clusters $k$ and $\ell$, respectively.
From the ILP formulation, we have the following constraints,
\begin{equation}\label{eq:nknl}
  n_k Q_{k\ell} = n_{\ell} Q_{\ell k}, \quad n_k \geq Q_{\ell k}, \quad n_{\ell} \geq Q_{k\ell}
\end{equation}
for each edge in the quotient graph.
Let $h = \text{gcd}(n_k,n_{\ell})$ so that $n_k = d_k h$ and $n_{\ell} = d_{\ell} h$, where $\text{gcd}(a,b)$ is the greatest common divisor of two integers $a$ and $b$.
Also, find $m$ such that,
\begin{equation}
  \frac{m}{c} = \frac{Q_{\ell k}}{n_k} = \frac{Q_{k \ell}}{n_{\ell}}
\end{equation}
Pick an ordered sequence of integers $b = \{b_1,b_2,\ldots,b_m\}$ such that,
\begin{equation}\label{eq:b}
  \sum_{j=1}^m b_j = h, \quad b_j \geq 1
\end{equation}
For each node $u_i \in \mathcal{C}_k$, $i = 0,\ldots,n_k-1$, we create edges $(u_i,w_{f(i,r_1,r_2)})$ for $r_1 = 0, \ldots, d_{\ell}-1$ and $r_2 = 1, \ldots, m$, where $w_j \in \mathcal{C}_{\ell}$, $j = 0, \ldots, n_{\ell}-1$.
\begin{equation}\label{eq:uw}
  \begin{aligned}
  f(i,r_1,r_2) = \left(i + r_1 c + \sum_{j=1}^{r_2} b_j \right) \bmod n_{\ell}
  \end{aligned}
\end{equation}
Note that each triplet $(i,r_1,r_2)$ yields a unique edge so that there are in total $n_k d_{\ell} m = n_k Q_{k\ell}$ edges, the required number of edges.
Also, each vertex $u_i \in \mathcal{C}_k$ is connected with $Q_{k\ell}$ vertices in $\mathcal{C}_{\ell}$.
Alternatively, we can generate the same set of edges with respect to the vertices in $\mathcal{C}_{\ell}$.
For each vertex $w_i \in \mathcal{C}_{\ell}$, $i = 0,\ldots,n_{\ell}-1$, we create edges $(w_i,u_{g(i,r_3,r_4)})$ for $r_3 = 0,\ldots,d_k-1$ and $r_4 = 1,\ldots,m$,
\begin{equation}\label{eq:wu}
  \begin{aligned}
    g(i,r_3,r_4) = \left( i + r_3 h + \sum_{j=1}^{r_4} b_{m-j+1} \right) \bmod n_k
  \end{aligned}
\end{equation}
Once again, note that each triplet ($i,r_3,r_4)$ yields a unique edge so that there are in total $n_{\ell} d_k m = n_{\ell} Q_{\ell k}$ edges, the same number as in Eq. \eqref{eq:uw} by the constraint in the ILP restated in Eq. \eqref{eq:nknl}.
Also, each vertex $w_i \in \mathcal{C}_{\ell}$ is connected with $Q_{\ell k}$ vertices in $\mathcal{C}_k$.
As both Eq. \eqref{eq:uw} and Eq. \eqref{eq:wu} create the same number of edges, we must  prove the statement that the two sets of edges are equal, that is, for every edge $(u_i,w_{f(i,r_1,r_2)})$ created by Eq. \eqref{eq:uw}, there is a corresponding edge $(w_{f(i,r_1,r_2)},u_i)$ created by Eq. \eqref{eq:wu}.
In other words, for every pair $(r_1,r_2)$, we must find a corresponding pair $(r_3,r_4)$ such that $g(f(i,r_1,r_2),r_3,r_4) = i$, or written with the definitions of $f$ and $g$,
\begin{equation}\label{eq:inter}
  \begin{aligned}
    \left[ \left( i + r_1 c + \sum_{j=1}^{r_2} b_j \right)\right. &\bmod n_{\ell} + r_3 c \\
    &\left.+ \sum_{j=1}^{r_4} b_{m-j+1} \right] \bmod n_k = i
  \end{aligned}
\end{equation}
Using the definition of the modulo, there exists some integer $q$ such that we may rewrite the term modulo $n_{\ell}$ as,
\begin{equation}\label{eq:modnl}
   \left( i + r_1 h + \sum_{j=1}^{r_2} b_j \right) \bmod n_{\ell} =  i + r_1 h + \sum_{j=1}^{r_2} b_j - qn_{\ell}
\end{equation}
There are two cases we must examine; if $r_2 < m$ and if $r_2 = m$.\\
\textbf{Case 1: $r_2 < m$.}
Set $r_4 = m - r_2$ so that the second summation becomes,
  \begin{equation}\label{eq:r2ltm}
    \sum_{j=1}^{m-r_2} b_{m-j+1} = \sum_{j = r_2+1}^m b_j
  \end{equation}
Thus the total summation becomes $\sum_{j=1}^{m} b_j = c$.
Substituting the results of Eqs. \eqref{eq:modnl} and \eqref{eq:r2ltm} into Eq. \eqref{eq:inter},
\begin{equation}
  \left[ i + (r_1 + r_3 + 1)h - q n_{\ell} \right] \bmod n_k
\end{equation}
From Eq. \eqref{eq:uw}, we know $0 \leq i \leq n_k-1$ so that $i \bmod n_k = i$, thus, we are left to find $r_3$ such that for some integer $t$,
\begin{equation}\label{eq:t1}
  \begin{aligned}
  (r_1 + r_3 + 1)h - q_{n\ell} &= tn_k\\
  &\Rightarrow \quad r_3 = td_k + q d_{\ell} - r_1 - 1
  \end{aligned}
\end{equation}
From Eq. \eqref{eq:wu}, we know $0 \leq r_3 \leq d_k-1$.
Applying the bounds and moving everything not dependent on $t$ to the expressions for the bounds,
\begin{equation}\label{eq:t1bnd}
  r_1 + 1 - qd_{\ell} \leq td_k \leq d_k + r_1 - qd_{\ell} 
\end{equation}
For there to certainly exist some integer $t$ that satisfies Eq. \eqref{eq:t1bnd}, the number of integers between the bounds, inclusive must be at least equal to $d_k$.
\begin{equation}
  (d_k + r_1 - qd_{\ell}) - (r_1 + 1 - qd_{\ell}) + 1 = d_k
\end{equation}
Thus there exists precisely one value of $t$ which satisfies Eq. \eqref{eq:t1}.\\
\textbf{Case 2: $r_2 = m$.}
For this case, set $r_4 = m$ so both summations are over the total sequence $b$, each of which totals $h$.
Using Eqs. \eqref{eq:modnl} and the fact that both summations equal $h$ allows us to rewrite Eq. \eqref{eq:inter} as,
\begin{equation}
  \left[ i + (r_1 + r_3 + 2)h - q n_{\ell}\right] \bmod n_k
\end{equation}
Once again, we must show that there exists some integer $t$ such that,
\begin{equation}
  \begin{aligned}
  (r_1+r_3 + 2)h - qn_{\ell} &= tn_k\\
  &\Rightarrow \quad r_3 = td_k + qd_{\ell} - r_1 - 2
  \end{aligned}
\end{equation}
Following the same procedure as in case 1, using the bounds $0 \leq r_3 \leq d_k - 1$, we see that $t$ must satisfy,
\begin{equation}
  r_1 + 2 - qd_{\ell} \leq td_k \leq d_k + r_1 + 1 - qd_{\ell}
\end{equation}
where the gap between the bounds, inclusive, is,
\begin{equation}
  (d_k + r_1 + 1 - qd_{\ell}) - (r_1 + 2 - qd_{\ell}) + 1 = d_k
\end{equation}
In summary, we have shown that Eqs. \eqref{eq:uw} and \eqref{eq:wu} each generate $n_kQ_{k\ell} = n_{\ell} Q_{\ell k}$ unique edges, and that for each edge in Eq. \eqref{eq:uw}, the same edge also appears in Eq. \eqref{eq:wu}, thus the sets of edges are equal.
\subsection{An Example}
A complete example of the process outlined in the previous sections \ref{ssec:feas}, \ref{ssec:ilp}, and \ref{ssec:wire} is shown in Fig. \ref{fig:example}.
The diagram of the quotient graph is shown in Fig. \ref{fig:example}(A) which consists of three vertices labeled $\mathcal{C}_1$, $\mathcal{C}_2$, and $\mathcal{C}_3$.
Vertices $\mathcal{C}_1$ and $\mathcal{C}_2$ have self-loops, $\mathcal{D}(\mathcal{C}_1) = 1$ and $\mathcal{D}(\mathcal{C}_2) = 2$.
There are two edges $\mathcal{F}_1 = (\mathcal{C}_1,\mathcal{C}_2)$ and $\mathcal{F}_2 = (\mathcal{C}_1,\mathcal{C}_3)$
The first edge has weights $\mathcal{W}(\mathcal{F}_1,0) = 2$ and $\mathcal{W}(\mathcal{F}_1,1) = 1$ and the second edge has weights $\mathcal{W}(\mathcal{F}_2,0) = 1$ and $\mathcal{W}(\mathcal{F}_2,1) = 2$.
The description of the quotient graph is summarized in the quotient adjacency matrix $Q$ also shown in Fig. \ref{fig:example}(A).
Additionally, the matrix $A$ as described in Section \ref{ssec:ilp} is shown which has a null space of dimension one.
The full ILP is shown in Fig. \ref{fig:example}(B) where the weights in the cost function, $\textbf{c} = (c_1,c_2,c_3)$, are chosen to all be one.
We choose to scale the solution by $s = 2$.
Note that as $\mathcal{D}(\mathcal{C}_1)$ is odd, that $2s x_1 = n_1$ while $sx_2 = n_2$ and $sx_3 = n_3$.
The lower bounds are found by using Eq. \eqref{eq:xlb}, $x_1^L = \frac{1}{2} \max \{2,2,1\}$, $x_2^L = \max \{3,2\}$, and $x_3^L = \max \{1,1\}$.
The solution to the ILP is, $x_1^* = 1$, $x_2^* = 4$, and $x_3^* = 1$, which can be converted to the cluster cardinalities $n_1 = 4$, $n_2 = 8$, and $n_3 = 2$.
In Fig. \ref{fig:example}(C) the resulting edges for each self-loop and edge in the original quotient graph are listed.
The edges created as prescribed by $\mathcal{D}(\mathcal{C}_1)$ are generated using Eq. \eqref{eq:intraodd} while the edges created as prescribed by $\mathcal{D}(\mathcal{C}_2)$ are generated using Eq. \eqref{eq:intraeven}.
Each edge created is sorted by the texture of its originating self-loop or edge in the original quotient graph.
Finally, a diagram of the resulting graph with $n=14$ vertices is shown in Fig. \ref{fig:example}(D) where the nodes are shaded and the edges are textured according to their originating feature in the quotient graph in Fig. \ref{fig:example}(A).
\subsection{Random Graphs with Non-Trivial MBC}\label{ssec:random}
In principle there may be more than one way to choose the sequence $b$ in Eq. \eqref{eq:b}.
After choosing one such sequence, the wiring procedure described in Eqs. \eqref{eq:intraeven}, \eqref{eq:intraodd}, \eqref{eq:uw}, and \eqref{eq:wu} is deterministic so for each quotient graph $\mathcal{Q}$ with cluster cardinalities $n_i$, $i = 1,\ldots,p$, the process will create one realization.
The graph created so far has the property that the partition of the nodes induced by the MBC $\mathcal{C}$ will be equal to the OAG, $\mathcal{O}$, due to the particular wiring of the edges discussed in section \ref{ssec:wire}.\\
\indent 
The procedure laid out in the previous subsections can be extended to the case that one is interested instead in generating a random graph with non-trivial MBCs (and not necessarily the automorphism group of the graph).
This can be done by randomly rewiring the edges of the network obtained in the procedure in sections \ref{ssec:feas}, \ref{ssec:ilp}, and \ref{ssec:wire} in such a way that the quotient graph is preserved, as described next.\\
\indent
For each set of intra-cluster edges in cluster $\mathcal{C}_k$, randomly choose $2$ edges, $(u_i,u_{i'})$ and $(u_j,u_{j'})$.
If $i \neq j'$ and $j \neq i'$, then remove these edges and add two new edges $(u_i,u_{j'})$ and $(u_j, u_{i'})$.
Repeat this process a suitable number of times.\\
\indent 
For each set of inter-cluster edges between cluster $\mathcal{C}_k$ and $\mathcal{C}_{\ell}$, randomly choose $2$ edges $(u_i,w_{i'})$ and $(u_j,w_{j'})$.
Remove these edges and add two new edges $(u_i,w_{j'})$ and $(u_j,w_{i'})$.
Repeat this process a suitable number of times.
\begin{algorithm}[H]
  \caption{Random Graph with Non-Trivial Balanced Coloring}
  \begin{algorithmic}[1]
    \REQUIRE $\mathcal{Q}$ is a quotient graph
    \REQUIRE $\textbf{c}$ is a positive vector of length $p$.
    \STATE Construct $A$ from $\mathcal{Q}$ according to Eq. \eqref{eq:Adef}
    \STATE Solve Eq. \eqref{eq:ilp} for $\textbf{x}$.
    \STATE If desired, scale $x_i \gets sx_i$ for integer $s > 1$ for $i=1,\ldots,p$
    \STATE $n_i = 2 x_i$ if $Q_{ii}$ is odd, and $n_i = x_i$ otherwise, $i = 1,\ldots,p$
    \FOR {$i = 1, \ldots, p$}
      \IF {$Q_{ii} > 0$}
        \STATE Use either Eq. \eqref{eq:intraeven} or Eq. \eqref{eq:intraodd} if $Q_{ii}$ is even or odd, respectively, to wire the intra-cluster edges.
        \STATE If desired, randomize the edges via swapping.
      \ENDIF
    \ENDFOR
    \FOR {All edges $\mathcal{F}_i = (\mathcal{C}_k,\mathcal{C}_\ell)$}
      \STATE Use Eq. \eqref{eq:uw} to wire the inter-cluster edges.
      \STATE If desired, randomize the edges via swapping.
    \ENDFOR
  \end{algorithmic}
\end{algorithm}
\section{Comparing the MBC and the OAG}
\begin{figure}
  \centering
  \includegraphics[width=\columnwidth]{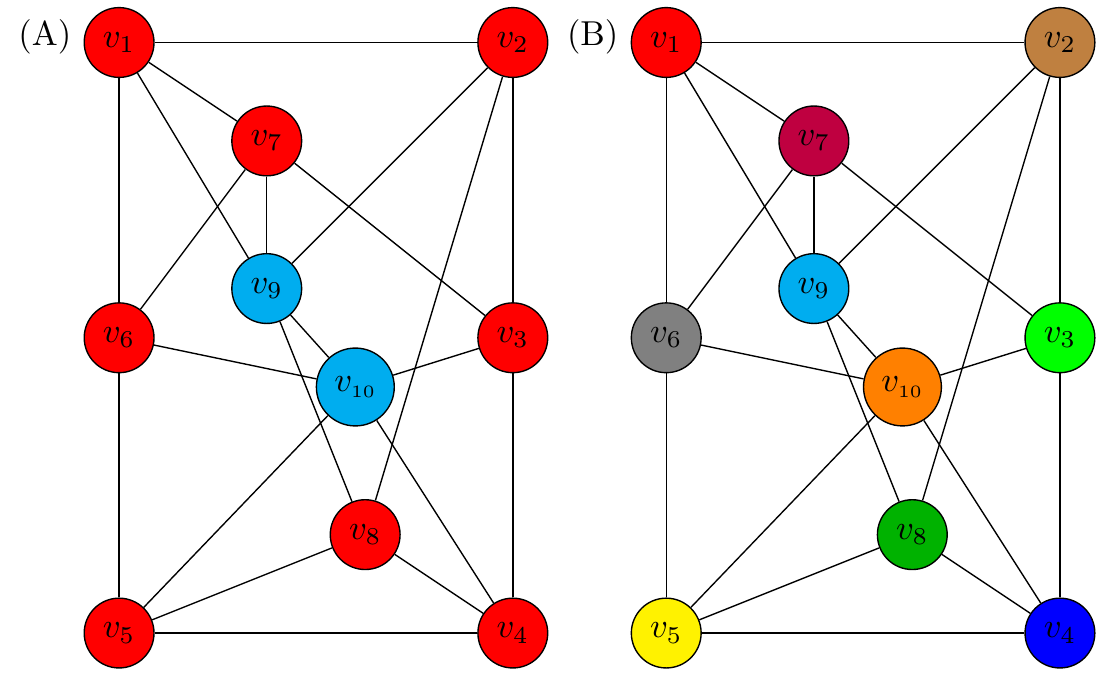}
  \caption{A network generated with two clusters in its MBC.
  In (A), the nodes are colored according to the MBC partition, where one can see that each red node is connected to three other red nodes and a blue node, while each blue node is connected to four red nodes and the other blue node.
  In (B), the nodes are colored according to the OAG partition.
  For this graph, each $\mathcal{O}_i \in \mathcal{O}$ consists of a single node.}
  \label{fig:trivial}
\end{figure}
\begin{figure}
  \centering
  \includegraphics[width=\columnwidth]{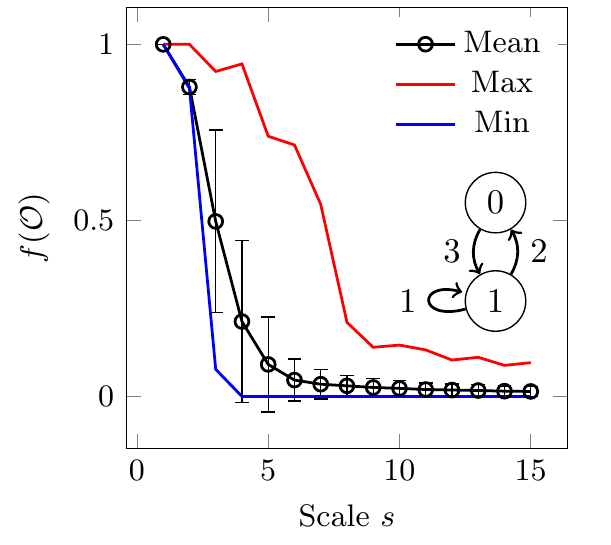}
  \caption{Comparing the cardinality of the OAG $\mathcal{O}$ to the cardinality of the MBC $\mathcal{C}$.
  Each graph is constructed using the quotient graph shown in the inset, so that $|\mathcal{C}| = 2$.
  The edge weights and self-loop weights are shown as well.
  The value of $s$ is the factor with which we scale the solution to Eq. \eqref{eq:ilp}.
  For each value of $s$ we generate 1000 graphs.
  The edges are randomly rewired preserving the MBC, and then the number of orbits of the automorphism group is determined.
  The black line is the mean with error bars representing one standard deviation.}
  \label{fig:size}
\end{figure}

It has been shown the OAG and MBC may not align \cite{siddique2018symmetry,schaub2016graph,kudose2009equitable}.
An example of this type of graph is presented in Fig. \ref{fig:trivial} where we show a network for which the equitable partition consists of two clusters, but the number of orbital partitions consists of ten clusters, the number of vertices, each consisting of a single vertex.\\
\indent 
We use the framework developed in this paper to numerically examine when to expect the MBC and the OAG to align and when they will not.
In all cases, we generate graphs $\mathcal{G}$ with the randomization procedure of section \ref{ssec:random} so that the MBC and OAG may or may not align.
We first examine how the size of the graph can affect $|\mathcal{O}|$ where we adjust the size by scaling the cluster cardinalities $\textbf{n}$ by a positive integer $s$.
We define the following metric,
\begin{equation}\label{eq:metric}
  f(\mathcal{O}) = \frac{n - |\mathcal{O}|}{n - |\mathcal{C}|}
\end{equation}
so that if $|\mathcal{O}| = |\mathcal{C}|$, $f(\mathcal{O}) = 1$, and if every orbit consists of a single vertex, $f(\mathcal{O}) = 0$.
Note that, by design, $|\mathcal{C}| < N$ so that, while the numerator may go to zero, the denominator does not change for a given quotient graph.
We choose a single quotient graph with only two vertices, $\mathcal{C} = \{\mathcal{C}_1, \mathcal{C}_2\}$, self-loops $\mathcal{D}(\mathcal{C}_1) = 0$ and $\mathcal{D}(\mathcal{C}_2) = 1$, and one edge $\mathcal{F} = \{\mathcal{F}_1\}$ where $\mathcal{F}_1 = (\mathcal{C}_1,\mathcal{C}_2)$ with weights $\mathcal{C}(\mathcal{F}_1,0) = 2$ and $(\mathcal{F}_1,1) = 3$, as shown in the inset of Fig. \ref{fig:size}.
After showing that this quotient graph is feasible and solving Eq. \eqref{eq:ilp} for the cluster cardinalities, $n_1 = 3$ and $n_2 = 2$, we are then free to scale $n_1$ and $n_2$ by any positive integer $s$ before wiring the graph.
For each value of $s$, we generate 1000 graphs, randomly rewiring the edges as described in section \ref{ssec:random}.

We see in Fig. \ref{fig:size} that if $s = 1$, i.e., each graph $\mathcal{G}$ is selected from the smallest graphs that can be represented by the quotient graph shown in the inset (following the procedure presented in section \ref{sec:results}), then the MBC and OAG almost always align.
As $s$ is increased, we see that $f(\mathcal{O})$ decreases rapidly, indicating that the MBC and OAG almost surely never align.
For $s > 6$, almost every graph generated has almost no non-trivial orbits of the automorphism group, that is $|\mathcal{O}| = n$.
While Fig. \ref{fig:size} shows results that are specific to the particular quotient graph in the inset, qualitatively similar behavior is seen for all quotient graphs examined.

\section{Conclusion}
Symmetries in complex networks and graphs have been shown to play an important role on the network dynamics, e.g., in the context of network synchronization \cite{siddique2018symmetry,schaub2016graph,nicosia2013remote,sorrentino2016complete}, and time averaged network dynamics \cite{siddique2017symmetries}.
However, to the best of our knowledge, no algorithms have been proposed that generate large networks with an assigned number of symmetries.
In this paper, we address this gap in the literature and propose a generating algorithm that is guaranteed to produce a network with an assigned number of symmetries from knowledge of a feasible quotient graph.
We also show how this algorithm can be extended to generate a graph with an assigned minimal balanced coloring (MBC), which may or may not coincide with the OAG of the graph.\\
\indent 
An analysis of anecdotal cases of networks has shown that the OAG and MBC of a graph may not always align \cite{siddique2018symmetry,kudose2009equitable}.
However, the question has remained unanswered of how common it is for the OAG and MBC of a graph to align.
Here we take advantage of our graph generating algorithm and show that when mapping a quotient network to a larger network with either a desired MBC or OAG, these two are never seen to align for large enough network size.
Our results indicate that the property that the MBC and the OAG of a graph may not align is indeed a generic feature of large graphs and networks.

\ifCLASSOPTIONcompsoc
  \section*{Acknowledgments}
\else
  \section*{Acknowledgment}
\fi

We would like to thank Lou Pecora, David Phillips, and Fabio Della Rossa for insightful conversations.
This work is supported by the National Science Foundation through NSF grant CMMI-1400193, NSF grant CRISP-1541148, and ONR Award No. N00014-16-1-2637 as well as HDTRA1-13-1-0020.

\ifCLASSOPTIONcaptionsoff
  \newpage
\fi

\bibliographystyle{IEEEtran}

\begin{IEEEbiography}
[{\includegraphics[width=1in,height=1.25in,clip,keepaspectratio]{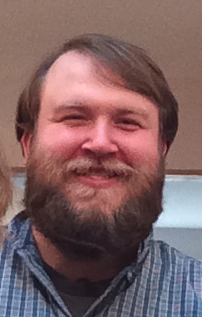}}]{Isaac Klickstein}
received his B.S. in mechanical engineering from the University of New Mexico in 2015.
He is currently a PhD student in the Department of Mechanical Engineering at the University of New Mexico.
His areas of specialty and interest are the control of complex networks, optimal control of nonlinear systems, and combinatorial problems on graphs.
\end{IEEEbiography}

\begin{IEEEbiography}
[{\includegraphics
[width=1in,height=1.25in,clip,keepaspectratio]{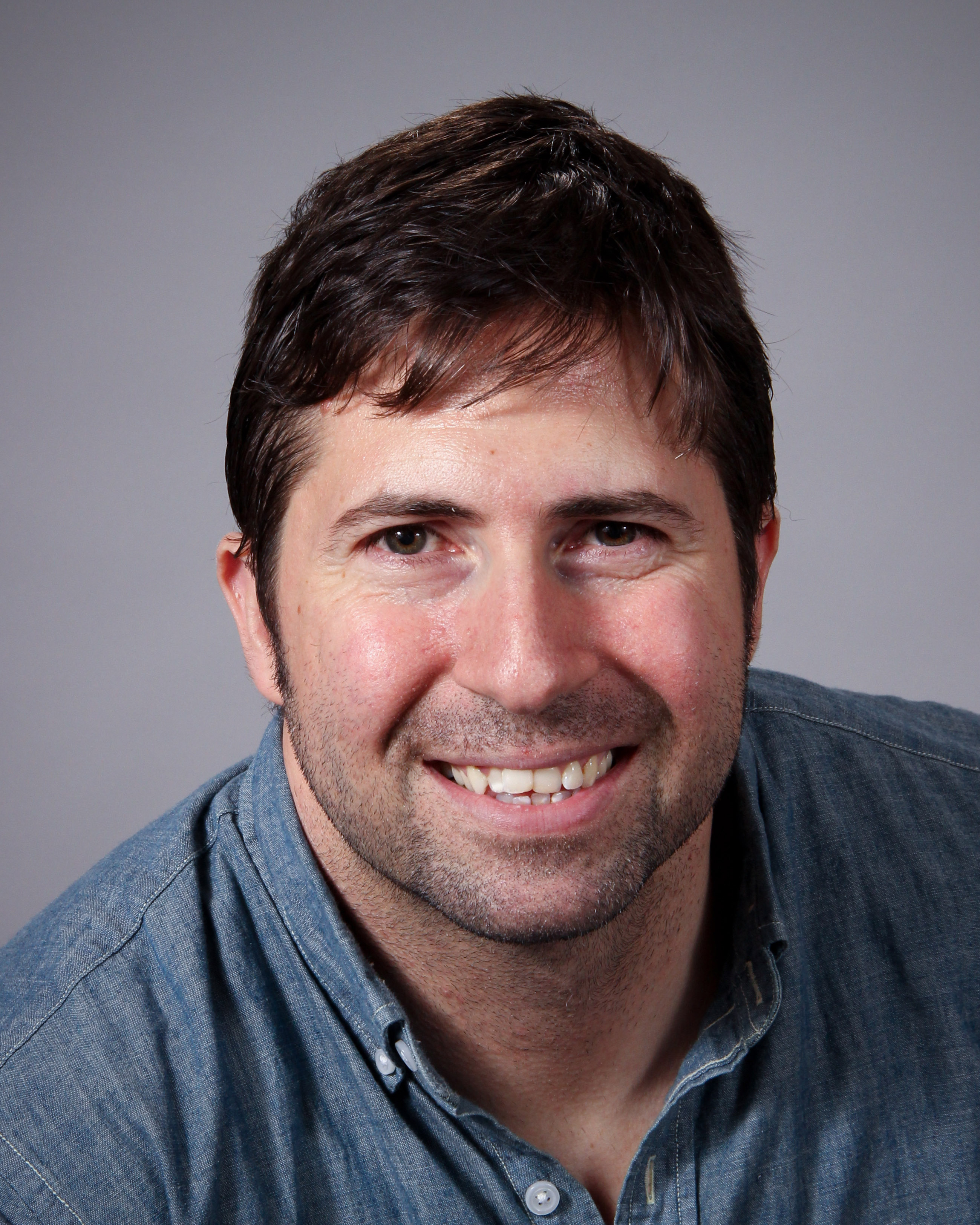}}]
{Francesco Sorrentino}
received a master's degree in Industrial Engineering from the University of Naples Federico II (Italy) in 2003 and a Ph.D. in Control Engineering from the University of Naples Federico II (Italy) in 2007.
His expertise is in dynamical systems and controls, with particular emphasis on nonlinear dynamics and adaptive decentralized control. His work includes studies on dynamics and control of complex dynamical networks and hypernetworks, adaptation in complex systems, sensor adaptive networks, coordinated autonomous vehicles operating in a dynamically changing environment, and identification of nonlinear systems. He is interested in applying the theory of dynamical systems to model, analyze, and control the dynamics of complex distributed energy systems, such as power networks and smart grids. Subjects of current investigation are evolutionary game theory on networks (evolutionary graph theory), the dynamics of large networks of coupled neurons, and the use of adaptive techniques for dynamical identification of communication delays between coupled mobile platforms.
He has published more than 40 papers in International Scientific Peer Reviewed Journals.
\end{IEEEbiography}

\end{document}